\documentclass[12pt]{amsart}
\usepackage[english]{babel}
\usepackage{egothic}
\usepackage[T1]{fontenc}

\usepackage{amsmath}
\usepackage{mathabx}
\usepackage{amsthm}
\usepackage{amssymb}
\usepackage[all]{xy}
\usepackage{mathrsfs}
\usepackage{enumerate}
\usepackage{physics}
\usepackage[notcite, final, notref]{showkeys}
\usepackage{dsfont}
\usepackage[dvips]{color}
\newtheorem{theorem}{Theorem}[section]

\newtheorem{lemma}[theorem]{Lemma}
\newtheorem{proposition}[theorem]{Proposition}
\newtheorem{corollary}[theorem]{Corollary}
\newtheorem{definition}[theorem]{Definition}

\theoremstyle{definition}
\newtheorem{remark}[theorem]{Remark}

\newcommand{\C}{{\mathbb C}}
\newcommand{\R}{{\mathbb R}}
\renewcommand{\i}{{\bf i}}

\newcommand{\e}{{\bf e}}

\newcommand{\vertiii}[1]{{\left\vert\kern-0.25ex\left\vert\kern-0.25ex\left\vert #1 
        \right\vert\kern-0.25ex\right\vert\kern-0.25ex\right\vert}}

\DeclareFontEncoding{FMS}{}{}
\DeclareFontSubstitution{FMS}{futm}{m}{n}
\DeclareFontEncoding{FMX}{}{}
\DeclareFontSubstitution{FMX}{futm}{m}{n}
\DeclareSymbolFont{fouriersymbols}{FMS}{futm}{m}{n}
\DeclareSymbolFont{fourierlargesymbols}{FMX}{futm}{m}{n}
\DeclareMathDelimiter{\VERT}{\mathord}{fouriersymbols}{152}{fourierlargesymbols}{147}

\newcommand{\BC}{\mathbb{B}\mathbb{C}}

\newcommand{\z}{\overline{z}}

\usepackage{xcolor}

\begin{document}
\title[Discrete Wiener algebra in the Bicomplex Setting]
{Discrete Wiener algebra in the bicomplex setting, spectral factorization with symmetry, and superoscillations}

\author[D. Alpay]{Daniel Alpay}
\address{(DA) Schmid College of Science and Technology \\
Chapman University\\
One University Drive
Orange, California 92866\\
USA}
\email{alpay@chapman.edu}

\author[I. Lewkowicz]{Izchak Lewkowicz}
\address{(IL) School of Electrical and Computer Engineering\\
Ben-Gurion University of the Negev\\ P.O.B. 653\\ Beer-Sheva, 84105\\
Israel}
\email{izchak@bgu.ac.il}

\author[M. Vajiac]{Mihaela Vajiac}
\address{(MV) Schmid College of Science and Technology \\
Chapman University\\
One University Drive
Orange, California 92866\\
USA}
\email{mbvajiac@chapman.edu}

\keywords{Bicomplex analysis}%
\subjclass{Primary 30G35; Secondary 47A57} %
\thanks{D. Alpay thanks the Foster G. and Mary McGaw Professorship in
  Mathematical Sciences, which supported his research}

\maketitle
\begin{abstract}
In this paper we present parallel theories on constructing Wiener algebras in the bicomplex setting. 
With the appropriate symmetry condition, the bicomplex matrix valued case can be seen as a complex valued case and, in this matrix valued case, we make the necessary connection between classical bicomplex analysis and complex analysis with symmetry.
We also write an application to superoscillations in this case.
\end{abstract}

\noindent {\em }
\date{today}





\section{Introduction}
\setcounter{equation}{0}
\subsection{Prologue}
In the present work we define and study the properties of the Wiener algebra in the framework of bicomplex numbers, i.e. $\mathbb{BC}$.
We here address $\mathbb{BC}$ questions in two ways: Each, on the expense of doubling the dimension, reduces the problem to the complex setting.
One is through the introduction of symmetry and the other goes along idempotents decomposition. In the end, we also write an application to superoscillations in this case.

To set the framework we first recall the definition of the matrix-valued Wiener algebra.
\begin{definition}
\label{p-cplx-Wiener}
The complex Wiener algebra $\mathcal W^{p\times p}$ consists of the matrix-valued functions of the form 
\begin{equation}
  \label{wiener-matrix}
  f(e^{it})=\sum_{n\in\mathbb Z} f_ne^{int}
\end{equation}
with $f_n\in\mathbb C^{p\times p}$ and
$\sum_{n\in\mathbb Z}\|f_n\|<\infty$ (where $\|f_n\|$ denotes the operator norm of the matrix $f_n$)
endowed with pointwise multiplication and norm
\begin{equation}
  \|f\|=\sum_{n\in\mathbb Z}\|f_n\|.
\end{equation}
\end{definition}

We note that in analysis quite often the Wiener algebra setting provides a convenient intermediate step, between the general
$\mathbf L_\infty$ case and the rational case (see Theorem \ref{Arche-de-la-defense} below for the latter). The algebra
$\mathcal W^{p\times p}$ is an important example of Banach algebra.
When $p=1$ this algebra is commutative, and will be denoted by $\mathcal W$.
Counterparts  and extensions of the Wiener algebra have been studied in various
other settings; see e.g. \cite{a-wiener,MR3554256,MR3688847}, \cite[\S II.1]{MR90c:47022}.\smallskip

 

In the scalar case (that is, $p=1$), the celebrated Wiener-L\'evy theorem asserts that point invertibility of an element of $\mathcal W$ implies invertibility in the algebra itself; see \cite[Lemma IIe p. 14]{wiener-1932}. We also note that in this case ($p=1$), the same theorem can be proved using Gelfand's theory of commutative Banach algebras and we write the statement as follows:

\begin{theorem}
Suppose that: $\displaystyle f(e^{i\theta})=\sum_{-\infty}^{\infty}c_n\, e^{in\theta}$,
with $\sum_{-\infty}^{\infty}|c_n|<\infty$, and $f(e^{i\theta})\neq 0$ for any $\theta\in\mathbb{R}$, then there exist $\{\gamma_n\}_{n\in\mathbb{Z}}$ such that:
\begin{equation*}
\frac{1}{f(e^{i\theta})}=\sum_{-\infty}^{\infty}\gamma_n\, e^{in\theta},
\end{equation*}
and $\sum_{-\infty}^{\infty}|\gamma_n|<\infty$.
\end{theorem}
We already note at this point that the above theorem yields a direct proof in the bicomplex case, which we present in Section~\ref{sec-4-2}.
The result still holds in the matrix-valued case by considering the determinant of the function. \smallskip

There exist two important subalgebras of $\mathcal W^{p\times p}$  used in spectral factorizations, namely $\mathcal W_+^{p\times p}$ and
$\mathcal W_-^{p\times p}$, defined as follows:

\begin{definition}
  The subalgebra $\mathcal W_+^{p\times p}$ consists of the elements $f\in\mathcal W^{p\times p}$ for which $f_n=0$ for $n<0$ in the representation
  \eqref{wiener-matrix}. Similarly,   the subalgebra $\mathcal W_-^{p\times p}$ consists of the elements $f\in\mathcal W^{p\times p}$ for which $f_n=0$
  for $n>0$ in the representation  \eqref{wiener-matrix}.
\end{definition}

A right (resp. left) spectral factorization of $f\in\mathcal W^{p\times p}$ is a representation
of the form $f=f_-f_+$ (resp. $f_+f_-$) where $f_+^{\pm 1}$ (resp. $f_-^{\pm 1}$) belongs to $\mathcal W_+^{p\times p}$ (resp.
$\mathcal W_-^{p\times p}$). We recall:

\begin{theorem}
\label{Wiener_gen}
Let $f\in\mathcal W^{p\times q}$  be a function with a non-identically vanishing determinant. If $f(e^{it})>0$ (pointwise, in the sense of positive matrices),
then it admits both right
and left spectral factorizations in $\mathcal W^{p\times q}$.
\end{theorem}

We note that, while $\mathcal W^{p\times p}$ is not commutative, factorization and inversion results are readily deduced from the scalar case using determinant; see e.g. \cite[\S 2]{DG1}. For completeness we also
mention the sources \cite{MR4:218g,MR0355675,ggk2} for more information on matrix-valued and operator-valued Wiener algebras. See in particular \cite[XXIX]{ggk2}.\smallskip


We define and study the Wiener algebra in the bicomplex ($\BC$) setting in two different ways, both reducing the problem to the complex setting.
One way, as in our previous paper \cite{ALV}, is through the introduction of a specific symmetry to complex matrices of twice the size,
which allows us to construct a bridge with the bicomplex setting. The other way 
goes along the idempotent decomposition of bicomplex numbers. In both cases, one doubles the dimension of the underlying framework, but in a different way.
The first approach, restricting a realization to be symmetric, comes to the expense of minimality, however we can write significant results in
this direction as well, as this second approach is classical in $\BC$ analysis. \smallskip

In this paper, our research focuses on three types of results in the theory of matrix Wiener algebra on the bicomplex space:\smallskip

$(i)$ Inversion Theorems,\\
$(ii)$ Factorization Theorems,\\
$(iii)$ Applications to Superoscillations.\smallskip

We present our results in Subsection~\ref{bcplx_alg_intro} and Subsection~\ref{outline_BC_matrix}, with complete proofs written in the subsequent parts of the paper. In order to do so, we first recall some properties of  $\BC$, the algebra of bicomplex numbers.

 
\subsection{ $\mathbb{BC}$ and the scalar bicomplex case:}
\label{bcplx_alg_intro}
We set the stage through the introduction of $\BC$, in the same fashion as~\cite{CSVV,bcbook,Price}, the key definitions and results for the case of holomorphic functions of a bicomplex variable, with a description of the bicomplex matrix setting following in Section~\ref{bcplx}.\smallskip

The algebra of bicomplex numbers, $\BC$, is generated by $i$ and $j$, two commuting imaginary units.  The product of the two commuting units  $i$ and $j$ is denoted by $ k := i\, j$ and we note that $k$ is a hyperbolic unit, i.e. it is a unit which squares to $1$.  Because of these various units in $\BC$, there are several
different conjugations that can be defined naturally. We will make use of these appropriate conjugations in this paper, and
we refer the reader to~\cite{bcbook,mltcplx} for more information on bicomplex and multicomplex analysis.\medskip

$\BC$ is not a division algebra: it has two distinguished zero divisors, $\e_1$ and $\e_2$,
\[
\e_1:=\frac{1+k}{2}\,,\qquad \e_2:=\frac{1-k}{2}\,,
\]
which are idempotent, linearly independent over the reals, and mutually annihilating with respect to the bicomplex multiplication:
\begin{align*}
    \e_1 \cdot \e_2 &= 0,\qquad
  \e_1^2=\e_1 , \qquad \e_2^2 =\e_2\,,\\
  \e_1 +\e_2 &=1, \qquad \e_1 -\e_2 = k\,.
\end{align*}
Just like $\{1, j \},$ they form a basis of the complex algebra
$\BC$, which is called the {\em idempotent basis}. If we define the
following complex variables in $\C(\i)$:
\begin{align}
  \label{l1l2}
  \lambda_1 := z_1-i z_2,\qquad \lambda_2 := z_1+i z_2\,,
\end{align}
the $\C(\i)$--{\em idempotent representation} for $Z=z_1+j z_2$ is
given by
\begin{align*}
  Z &= \lambda_1\e_1+\lambda_2\e_2\,.
\end{align*}

Simple algebra yields:
\begin{equation}
\begin{split}
  z_1&=\frac{\lambda_1+\lambda_2}{2}\\
    z_2&=\frac{i(\lambda_1-\lambda_2)}{2}.
  \end{split}
  \end{equation}

The $\C(\i)$--idempotent representation is the only one for which
multiplication is component-wise, as shown in the next lemma.

\begin{lemma}
  \label{prop:idempotent}
  The addition and multiplication of bicomplex numbers can be realized
  component-wise in the idempotent representation above. Specifically,
  if $Z= \lambda_1\,\e_2 + \lambda_2\,\e_2$ and $W= \mu_1\,\e_1 + \mu_2\,\e_2 $ are two
  bicomplex numbers, where $a_1,a_2,b_1,b_2\in\C(\i)$, then
  \begin{eqnarray*}
    Z+W &=& ( \lambda_1+\mu_1)\,\e_1  + ( \lambda_2+\mu_2)\,\e_2   ,  \\
    Z\cdot W &=& ( \lambda_1\mu_1)\,\e_1  + ( \lambda_2 \mu_2)\,\e_2   ,  \\
    Z^n &=&  \lambda_1^n \,\e_1  +  \lambda_2^n \,\e_2  .
  \end{eqnarray*}
  Moreover, the inverse of an invertible bicomplex number
  $Z= \lambda_1\e_1 +  \lambda_2\e_2 $ (in this case $ \lambda_1 \cdot  \lambda_2 \neq 0$) is given
  by
  $$
  Z^{-1}=  \lambda_1^{-1}\e_1 +  \lambda_2^{-1}\,\e_2 ,
  $$
  where $ \lambda_1^{-1}$ and $ \lambda_2^{-1}$ are the complex multiplicative
  inverses of $ \lambda_1$ and $ \lambda_2$, respectively.
\end{lemma}

One can see this also by computing directly which product on the
bicomplex numbers of the form
\begin{align*}
  x_1 + i x_2 + j x_3 + k x_4,\qquad x_1,x_2,x_3,x_4\in\R
\end{align*}
is component wise, and one finds that the only one with this property
is given by the mapping:
\begin{align}
  \label{shakira}
  x_1 + i x_2 + j x_3 + k x_4 \mapsto ((x_1 + x_4) + i (x_2-x_3), (x_1-x_4) + i (x_2+x_3))\,,
\end{align}
which corresponds to the idempotent decomposition
\begin{align*}
  Z = z_1 + j z_2 = (z_1-i z_2)\e_1 + (z_1+i z_2)\e_2\,,
\end{align*}
where $z_1 = x_1+i x_2$ and $z_2 = x_3+i x_4$.\\

A special subset of the bicomplex space is defined here:
\begin{definition}
  The space of {\em hyperbolic numbers $\mathbb D$} is defined as the subset of the bicomplex space $\BC$ of the form $Z=a+bk$, with $a$ and $b$ in $\mathbb R$. 
  \end{definition}

\begin{definition}
The bicomplex number $Z$ will be said to be positive if both $\lambda_1$ and $\lambda_2$ in~\eqref{l1l2} are positive real numbers in the idempotent representation, and we write $Z\in\mathbb{D}_+$. The space $\mathbb{D}_+$ is called the cone of hyperbolic positive numbers.
  \end{definition}

Using the idempotent representation one can define three norms associated with it, see ~\cite{SVV}. 
We will make use of all of these norms in our work, according to the types of properties that we wish to analyze.
The first two definitions represent the Lie norm and the dual Lie norm as follows:
\begin{definition}
\label{dual_Lie}
The dual Lie norm is, up to a factor of $2$:
$$||Z||=|\lambda_1|+|\lambda_2|=2\mathcal{L}^*(Z),$$ 
\end{definition}
and
\begin{definition}
The Lie norm is:
$$\mathcal{L}(Z)=\max\{|\lambda_1|,|\lambda_2|\}.$$ 
\end{definition}

A slightly non-standard norm is:

\begin{definition}
The hyperbolic valued norm:
$\|Z\|_{\mathbb{D}_+}=|\lambda_1|\e_1+|\lambda_2|\e_2$. 
\end{definition}
It is easy to see that $\|Z\|_{\mathbb{D}_+}\in\mathbb{D}_+$ .
We will need the following results regarding the various norms:
\begin{proposition}
The norm in Definition~\ref{dual_Lie} is sub-multiplicative, i.e.
\begin{equation}
  ||ZW||\le ||Z||\cdot ||W||.
\end{equation}
\end{proposition}
\begin{proof}
  \[
    Z=\lambda_1e_1+\lambda_2 e_2\quad {\rm and}  \quad   Z=\mu_1e_1+\mu_2 e_2.
  \]
  Then, $ZW=\lambda_1\mu_1e_1+\lambda_2\mu_2e_2$ and
  \[
||ZW||=|\lambda_1\mu_1|+|\lambda_2\mu_2|\le (|\lambda_1|+|\lambda_2|)(|\mu_1|+|\mu_2|)=||Z||\cdot ||W ||.
    \]
  
  \end{proof}

For $\mathbb{BC}$-valued functions, the natural counterpart of the bi-disk in $\mathbb C^2$ is:
\begin{equation}
  \label{bc-bi-disk}
  \mathbb K=\left\{(z_1,z_2)\in\mathbb C^2\,;\, |z_1\pm iz_2|<1\right\}.
\end{equation}

We note that $\mathbb K$ contains $B(0,1)\times \left\{0\right\}$,and so the single variable theory is obtained as a special case.

Here we discuss the boundaries and properties of this bicomplex domain $\mathbb K$.
\begin{lemma}
\label{bcpx_trig}
  The distinguished boundary $|\lambda_1|=|\lambda_2|=1$ of the bidisk in the decomposition along the idempotents corresponds to
  \begin{equation}
    \label{ab-ab}
    \left\{ e^{it}(\cos s, \sin s)\, ;t,s\in[0,2\pi]\right\}, 
  \end{equation}
  i.e.
  \begin{equation}
    \label{z-b-a-1}
    Z=e^{it}e^{js}
    \end{equation}
\end{lemma}
\begin{proof}
  Since
  \[
|z_1\pm iz_2|^2=|z_1|^2+|z_2^2\mp 2{\rm Re}\,(iz_1\overline{z_2})|
    \]
the conditions $|z_1\pm iz_2|=1$ are equivalent to the two conditions
  \[
    \begin{split}
      |z_1|^2+|z_2|^2&=1\\
      {\rm Im}\, (z_1\overline{z_2})&=0.
    \end{split}
  \]
  Setting $z_1=uz_2$ when $z_2\not=0$, $u\in\mathbb R$, we get
\[
\partial\mathbb K=  \left\{(e^{i\theta},0),
\theta\in[0,2\pi]\right\}\bigcup\left\{\bigcup_{u\in\mathbb R}\left\{\left(\frac{ue^{i\theta}}{\sqrt{u^2+1}},\frac{e^{i\theta}}{\sqrt{u^2+1}}\right),\,
    \theta\in[0,2\pi]\right\}\right\},
\]
which is rewritten as \eqref{ab-ab}. To conclude, note that $\partial \mathbb K$ 
\[
z_1+jz_2=e^{it}(\cos s +j\sin s)=e^{it}e^{js}.
\]  
\end{proof}

We will call \eqref{ab-ab} the distinguished boundary of $\mathbb K$ and denote it by $\partial\mathbb K$.

The counterpart of the unit circle can now be defined:

\begin{definition}
The bicomplex unit circle $\partial \mathbb K$ is:
\begin{equation}
    \label{z-b-a}
  \partial\mathbb K=\left\{e^{it}e^{js},\, t,s\in[0,2\pi]\right\}
\end{equation}
 \end{definition}
The following two lemmas show us the structure of $\partial \mathbb K$:
 \begin{lemma}
 In the idempotent representation, the bicomplex unit circle becomes:
 \begin{equation}
  \partial\mathbb K=\left\{e^{i(t+s)}\e_1 + e^{i(t-s)}\e_2,\, t,s\in[0,2\pi]\right\}
\end{equation}
 \end{lemma}
 \begin{proof}
 Since:
 \[
  z_1=e^{it}\cos s\quad{\rm and}\quad z_2=e^{it}\sin s,
\]
we have that $\lambda_1$ and $\lambda_2$, defined by \eqref{l1l2}, are equal to
    \begin{equation}
      \label{dupont-circle-1}
    \lambda_1=e^{i(t-s)}\quad{\rm and}\quad     \lambda_2=e^{i(t+s)}
  \end{equation}
  and
  \[
    z_1+jz_2=e^{i(t-s)}e_1+e^{i(t+s)}e_2
  \]
 
 \end{proof}
\begin{lemma}
Elements of $\partial \mathbb K$ are invertible in $\mathbb{BC}$.
\end{lemma}

\begin{proof}
This follows from \eqref{z-b-a}.
\end{proof}

\begin{definition}
  Let the bicomplex Wiener algebra, denoted $\mathfrak W$  be the set of functions $f:\partial\mathbb K \to\mathbb{BC}$, such that:
  \begin{equation}
    \label{place-de-la-Bastille}
    f(Z)=\sum_{n\in \mathbb Z}f_nZ^n,\quad Z\in\partial\mathbb K,
  \end{equation}
  such that $\sum_{n\in\mathbb Z} ||f_n||<\infty$.
\end{definition}

We can also define the counterparts of  $\mathfrak W_{+}$ and   $\mathfrak W_{ -}$ in the usual way:
\begin{definition}
Let $\mathfrak W_+$  be the set of functions $f:\partial\mathbb K \to\mathbb{BC}$, such that:
  \begin{equation}
    f(Z)=\sum_{n\ge 0}f_nZ^n,\quad Z\in\partial\mathbb K,
  \end{equation}
  such that $\sum_{n\ge 0} ||f_n||<\infty$.
  \
  
  Then, let $\mathfrak W_-$  be the set of functions $f:\partial\mathbb K \to\mathbb{BC}$, such that:
  \begin{equation}
    f(Z)=\sum_{n\ge 0}f_{-n}Z^{-n},\quad Z\in\partial\mathbb K,
  \end{equation}
  such that $\sum_{n\ge 0} ||f_{-n}||<\infty$.
\end{definition}


In the case of the scalar bicomplex valued case in Subsection~\ref{scalar_W_BC} we prove:

\begin{theorem}
Assume that $f(Z)$ invertible on $\partial\mathbb K$. Then $f$ is invertible in the bicomplex Wiener algebra.
\end{theorem}

\begin{theorem}
Assume that $f(Z)>0$ invertible on $\partial\mathbb K$. Then $f$ admits left and right spectral factorizations.
  \end{theorem}

\subsection{Outline of the matrix valued case}
 \label{outline_BC_matrix}

\
In Section~\ref{Wiener_BC_matrix} we define the Wiener algebra in the bicomplex matrix valued case as:

 \begin{definition}
The bicomplex Wiener algebra $\mathfrak W^{p\times p}$ consists of the matrix-valued functions of the form 
\begin{equation}
  \label{bc_wiener-matrix}
  f(Z)=\sum_{n\in\mathbb Z} f_nZ^n
\end{equation}
with $f_n\in\BC^{p\times p}$ and
$\sum_{n\in\mathbb Z}||f_n||<\infty$.

Here $||f_n||$ denotes the following norm of the bicomplex matrix $f_n=f_{n,1}\e_1+f_{n,2}\e_2$:
$||f_n||=\|f_{n,1}\|+\|f_{n,2}\|$, where $\|f_{n,1}\|$ and $\|f_{n,2}\|$ are the operator norms as complex matrices.

This Wiener algebra is endowed with pointwise multiplication and norm
\begin{equation}
  \|f\|=\sum_{n\in\mathbb Z}\|f_n\|.
\end{equation}
\end{definition}

The elements of $\mathfrak W^{p\times p}_+$  are  of the form 
\begin{equation}
  \label{bc_wiener-matrix-+}
  f(Z)=\sum_{n\ge 0} f_nZ^n,
\end{equation}
while elements of $\mathfrak W^{p\times p}_-$  are of the form
 \begin{equation}
  \label{bc_wiener-matrix_-}
  f(Z)=\sum_{n\ge 0} f_{-n}Z^{-n}.
\end{equation}

We prove the following invertibility and factorization theorems (Theorem~\ref{p-bcplx-wiener} and ~\ref{bc_m_sym_complex}):
\begin{theorem}
Let $f\in \mathfrak W^{p\times p}$ be such that $f(Z)$ is invertible for $Z\in\partial \mathbb K$. Then $f$ is invertible in
  $\mathfrak W^{p\times p}$.
\end{theorem}

\begin{theorem}
  Let $f\in \mathfrak{W}^{p\times p}$ be such that $f(Z)>0$ for $Z\in\partial \mathbb K$. Then
  $f$ admits a spectral factorization $f(Z)=f_+(Z)f_+(Z)^*$ where $f_+^{\pm1}\in \mathfrak W^{p\times p}_+$.
\end{theorem}

In Subsection~\ref{bcplx_wiener_sharp} we also prove an equivalence between this bicomplex Wiener algebra and one induced by a special type of symmetry, called $\sharp$ symmetry in the complex matrix domain. We will have the same type of factorization result in this case, see Theorem~\ref{wiener_m_sym_bicomplex-2}.
For an equivalence theorem between the two algebras see Theorem~\ref{W_BC_equiv}.
In the same Section, in Subsection~\ref{rational_bcplx_sharp}, we address the rational bicomplex case with $\sharp$ symmetry as well.

\
Last but not least, as an applications of these results, in Section~\ref{Bcplx_sup} we introduce a theory of superoscillations in the bicomplex case.


\section{Matrix-Valued Complex Wiener Algebras}
\setcounter{equation}{0}
 

\subsection{Wiener Algebra and Matrix Symmetries}
\label{Wiener_sym}
\

We first describe relations in the Wiener algebras induced by symmetries with respect to invertible matrices. 
\begin{definition}
Let $Y\in GL_n(\mathbb C)$ and we define
\begin{equation}
M_Y=YMY^{-1},\quad M\in\mathbb C^{n\times n}.
  \end{equation}
  \end{definition}
  \begin{remark}
  It is immediate that, for $M,N\in\mathbb C^{n\times n}$ and $\lambda\in\mathbb C$,
  \begin{eqnarray*}
    (MN)_Y&=&M_YN_Y\\
    (M+N)_Y&=&M_Y+N_Y\\
    (\lambda I_n)_Y&=&\lambda I_n.
              \end{eqnarray*}
  \end{remark}
  We now have:
  \begin{lemma}
              Let now $f\in\mathcal W^{p\times p}$ be of the form \eqref{wiener-matrix}. if we define $\displaystyle f_Y(e^{it})=\sum_{n\in\mathbb Z}(f_n)_Ye^{int}$
              it follows that
\begin{equation*}
f_Y(e^{it})=(f(e^{it}))_Y.
  \end{equation*}
  \end{lemma}
From the uniqueness of the Fourier coefficients and Theorem~\ref{Wiener_gen} the following result follows:
\begin{proposition}
    \label{proper-Y}
  Let $f,g\in\mathcal W^{p\times p}$. We have:
    \begin{eqnarray*}
      (fg)_Y(e^{it})&=&f_Y(e^{it})g_Y(e^{it})\\
      f_Y(e^{it})&=&(f(e^{it}))_Y\\
      ((f(e^{it}))_Y)^*&=&((f(e^{it}))^*)_Y.
    \end{eqnarray*}
    Furthermore, $f\in\mathcal W_+^{p\times p}$ if and only if $f_Y\in\mathcal W_+^{p\times p}$.
    \end{proposition}

 \begin{definition}
The function $f\in\mathcal  W^{p\times p}$ is called $Y$-symmetric if $f(e^{it})=f_Y(e^{it})$, $t\in[0,2\pi]$.
\end{definition}

From uniqueness of the Fourier coefficients we obtain:

\begin{proposition}
  The function $f\in\mathcal  W^{p\times p}$ is $Y$-symmetric if and only if its Fourier coefficients are $Y$-symmetric:
  \[
f_n=(f_n)_Y,\quad n\in\mathbb Z.
    \]
  \end{proposition}

\begin{theorem}
\label{wiener_m_sym_complex}
  Let $f\in\mathcal W^{p\times p}$ be $Y$-symmetric for some $Y\in\mathbb{C}^{p \times p}$ and such that $f(e^{it})>0$ for $t\in[0,2\pi]$. Then, there is a uniquely determined spectral
  factorization $f=f_+f_+^*$ where $f_+\in\mathcal W_+^{p\times p}$ is $Y$-symmetric and satisfies $f_+(1)=I_{p}$ and invertible and any  other spectral factorization differs by an $Y$-symmetric
  unitary multiplicative constant on the right from $f_+$.
\end{theorem}
\begin{proof}
We leave the proof to the reader, since it is very similar to the proof of theorem~\ref{wiener_m_sym_complex-2} and uses parallel techniques.
\end{proof}
This above $Y$ symmetry does not apply for the more general, non-square case. We introduce another symmetry that will overcome this limitation.

\subsection{Wiener algebra with special symmetry}
\

We now turn towards a type of symmetry that will be useful in the bicomplex case and could be applied to the non-square case as well. 

Let $J\in \mathbb C^{2\times 2}$ be a matrix such that $J\, J^*=I_2$.
\begin{definition}
\label{bc_M_sharp}  
Let $M\in \mathbb C^{2\times 2}$. 
We define the $\sharp$-conjugate with respect to $J$ to be:
\[
M^\sharp=JMJ^*
\]
More generally, given $M\in\mathbb C^{2a\times 2b}$ a block matrix, we define
\begin{equation}
M^\sharp=J_aMJ_b^*,
\end{equation}
with
\begin{equation}
  \label{J_a=}
 J_a=J \otimes I_a,\qquad J_a^*=J^* \otimes I_a
\end{equation}
and where $\otimes$ denotes the  Kronecker (tensor) product of matrices.
\end{definition}

\begin{proposition}
  Let $M\in\mathbb C^{2p\times 2q}$ and $N\in\mathbb C^{2q\times 2r}$ be matrices of possibly different sizes
  Then
  \[
  (MN)^\sharp=M^\sharp N^\sharp.
  \]
    Let $M,N\in\mathbb C^{2p\times 2q}$. Then,
\[    
    (M+N)^\sharp=M^\sharp+N^\sharp.\\
\]
  \end{proposition}
\begin{definition}
For $f\in\mathcal W^{2p\times 2p}$ of the form \eqref{wiener-matrix}, we define its $\sharp-$ symmetric function to be:
              \begin{equation}
f^\sharp(e^{it})=\sum_{n\in\mathbb Z}(f_n)^\sharp e^{int}.
\end{equation}
\end{definition}
\begin{remark}
It is easy to see that:
\begin{equation}
f^\sharp(e^{it})=(f(e^{it}))^\sharp,
  \end{equation}
  therefore $f^{\sharp}\in \mathcal W^{2p\times 2p}$ of the form \eqref{wiener-matrix}.
\end{remark}

The counterpart of Proposition \ref{proper-Y} is:

  \begin{proposition}
  \label{proper-X}
    Let $p,q,s\in\mathbb N$ and let $f\in\mathcal W^{2p\times 2q}$ and $g\in\mathcal W^{2q\times 2s}$. We have:
    \begin{eqnarray}
      (fg)^\sharp(e^{it})&=&f^\sharp(e^{it})g^\sharp(e^{it})\\
      f^\sharp(e^{it})&=&(f(e^{it}))^\sharp\\
      ((f(e^{it}))^\sharp)^*&=&((f(e^{it}))^*)^\sharp.
    \end{eqnarray}
    Furthermore, 
$f\in\mathcal W_+^{2p\times 2q}$ if and only if $f^\sharp\in\mathcal W_+^{2p\times 2q}$.
    \end{proposition}


\begin{definition}
The matrix $M\in\mathbb C^{2\times 2}$ is called $\sharp$-symmetric with respect to $J$ if and only if it satisfies $M=M^\sharp$.
More generally the block matrix $M\in\mathbb C^{2a\times 2b}$ is called $\sharp$-symmetric with respect to $J$ if $M=M^\sharp$. This will hold if and only if all its $\mathbb C^{2\times 2}$ block entries
are $\sharp$-symmetric with respect to $J$.
\end{definition}

\begin{definition}
The function $f\in\mathcal  W^{2p\times 2p}$ is called $\sharp$-symmetric if $f(e^{it})=f^\sharp(e^{it})$, $t\in[0,2\pi]$.
\end{definition}

From the uniqueness of the Fourier coefficients follows directly the following result:

\begin{proposition}
  The function $f\in\mathcal  W^{2p\times 2p}$ is $\sharp$-symmetric if and only if its Fourier coefficients are $\sharp$-symmetric:
  \[
f_n=(f_n)^\sharp,\quad n\in\mathbb Z.
    \]
  \end{proposition}

\begin{theorem}
\label{wiener_m_sym_complex-2}
  Let $f\in\mathcal W^{2p\times 2p}$ be $\sharp$-symmetric and such that $f(e^{it})>0$ for $t\in[0,2\pi]$. Then, there is a uniquely determined spectral
  factorization $f=f_+f_+^*$ where $f_+\in\mathcal W_+^{2p\times 2p}$ is $\sharp$-symmetric and satisfies $f_+(1)=I_{p}$ and invertible and any  other spectral factorization differs by an $\sharp$-symmetric
  unitary multiplicative constant on the right from $f_+$.
\end{theorem}

\begin{proof}
  Using Proposition \ref{proper-Y} we can write:
  \[
      f(e^{it})=f_+(e^{it})f_+(e^{it})^*
    \]
    and
    \[
      \begin{split}
        f(e^{it})&=(f(e^{it}))^\sharp\\
        & =       (f_+(e^{it})f_+(e^{it})^*)\\
        &=(f_+(e^{it}))^\sharp(f_+(e^{it})^*)^\sharp\\
        &=(f_+(e^{it}))^\sharp((f_+(e^{it}))^\sharp)^*
        \end{split}
      \]
      By uniqueness up to left multiplicative factor of the spectral factor we have
      \[
        f_+(e^{it})=(f_+(e^{it}))^\sharp\cdot U
      \]
      and so
      \begin{equation}
        \label{yuiop}
f_+(z)=(f_+(z))^\sharp\cdot U, \quad |z|<1.
\end{equation}
for some $U\in GL_p(\mathbb C)$ which is moreover unitary.  To show that $U=I_n$ remark that $f_+(0)$ is invertible and so, with
$f_+(z)=\sum_{n=0}^\infty a_nz^n$ we have $a_0\in GL_p(\mathbb C)$. Setting $z=0$ in \eqref{yuiop} leads to $a_0=a_0U$ and so
$U=I_p$ since $a_0$ is invertible.
        \end{proof}

In the case the bicomplex algebra, as seen in~\cite{ALV} we have an important structural matrix:
\begin{equation}
\label{bc-J}
J=\begin{pmatrix} 0 &-1\\ 1&0\end{pmatrix}\, ,
\end{equation}
which has the obvious property that $J\, J^*=I_2$.

\begin{definition}
\label{bcplx-sharp-2by2}
The matrix $M\in\mathbb C^{2\times 2}$ is called bicomplex $\sharp$-symmetric if it satisfies $M=M^\sharp$, with respect to $J$ defined in ~\eqref{bc-J} above.
\end{definition}

As in~\cite{ALV} we have:
\begin{remark}
For  $J$ as in~\eqref{bc-J}, we have that $M$ is $J$-symmetric if and only if it is of the form
\begin{equation}
  \label{petunia}
M=\begin{pmatrix}z_1&-z_2\\ z_2&z_1\end{pmatrix}\,;\, z_1,z_2\in\mathbb C ,
\end{equation}
which, in Section~\ref{bcplx_wiener_sharp} will be discussed at length in relation with bicomplex analysis.
\end{remark}

\begin{remark}{\rm        
As we will see in Section~\ref{bcplx_matrix_sharp}, a bicomplex number can be viewed as a matrix of the form \eqref{petunia} (this gives in fact an explicit construction of the bicomplex
numbers in terms of complex matrices).
Bicomplex matrices will be identified with $\sharp$-symmetric complex ones of the appropriate dimension. This allow us follow the strategy introduced in our previous
work  \cite{ALV}, where we reduce the various results to ones in a complex Wiener algebra with symmetry.
Note that we will also use the second approach, splitting the problem into two classical complex ones via the idempotent representation. 
}
\end{remark}

    
\subsection{Rational Functions and the Wiener Algebra}
\label{Rational_Wiener}
\
In the general theory of  rational matrix-valued functions, we can speak about several types of realizations. We recall results that pertains to our work in the bicomplex case.

\begin{theorem}
  \label{real}
Every matrix-valued (say $\mathbb C^{a\times b}$-valued) rational function which is analytic at infinity admits a realization centered at infinity,
i.e. can be written in the form:
\[
  f(z)=D+C(zI_N-A)^{-1}B.
\]
\end{theorem}

\begin{remark}
If $f$ is not analytic at infinity, then it has a polynomial part which will not be of this form, i.e. there exists $P$ polynomial such that:
\[
  f(z)=P(z)+ D+C(zI_N-A)^{-1}B.
\]
\end{remark}

The realization is called minimal if $N$ is minimal.
We recall a proof of the existence of a (possibly non minimal) realization. The arguments will be used in the sequel to prove
the main result of the section.

\begin{proof}[Proof of Theorem \ref{real}]
Let $p_1,\ldots, p_M$ be the poles of $f$ (in the most elementary sense, meaning points where any entry of $f$
has a non-removable singularity). At each of these points $f$ has a finite Laurent expansion and we can write (since $f$ is assumed analytic at infinity)
  \[
f(z)=D+\sum_{m=1}^M\underbrace{\left(\sum_{k=1}^{k_m}\frac{H_{k,m}}{(z-p_m)^k}\right)}_{f_m(z)},\quad z\not\in\left\{p_1,\ldots, p_M\right\},
\]
where $D$ and the coefficients $H_{k,m}$ belong to $\mathbb C^{a\times b}$. We now write
  \begin{equation}
    f_m(z)=C_m(z I_{N_m}-A_m)^{-1}B_m
  \end{equation}
  where
  \begin{equation}
C_m=\underbrace{\begin{pmatrix}I_{a}&0_{a}&\cdots&0_a\end{pmatrix}}_{k_m\, a\times a \,\mbox{{\rm blocks}}} 
\end{equation}
\begin{equation}
  B_m=\begin{pmatrix}H_1\\ \vdots\\ H_{k_m}\end{pmatrix}
\end{equation}

and
\begin{equation}
  A_m=\begin{pmatrix}p_m I_a&I_a&0_a&\cdots&0\\
    0_a&p_m I_a&I_a&\cdots&0_a\\
    & & & & \\
    0_a&\cdots&&  p_m I_a&I_a\\
    0_a&\cdots&& 0_a & p_m I_a
      \end{pmatrix},
\end{equation}
where $A_m\in\C^{k_m a\times k_m a}$.
\
Thus  neither the realization of $f_m$ nor the realization of  $f$ (given by $f(z)=D+C(I_N-A)^{-1}B$) are necessarily minimal, where 
\begin{equation}
  C=\begin{pmatrix}C_1&\cdots &C_{M}\end{pmatrix},\quad A={\rm diag}(A_1,\ldots, A_{M}),\quad
  B=\begin{pmatrix}B_1\\ \vdots \\B_{M}\end{pmatrix}.
\end{equation}
\end{proof}

We recall the following result. See \cite{gk1} for more details and information. Note that we do not require the realization to be minimal. The proof
in \cite{gk1} still holds since $A$ is assumed to have to spectrum on the unit circle. The realization is then called {\sl regular} on the unit circle.

\begin{theorem}
\label{Riesz_proj}
  Let $f$ be a $\mathbb C^{p\times p}$-valued rational function. If $f$ has no singularities on the unit circle, then $f\in\mathcal W^{p\times p}$.
  Let $f(z)=D+C(zI_N-A)^{-1}B$ be a realization of $f$, possibly not minimal but with no spectrum on the unit circle, and let $P$ denote the Riesz projection corresponding to the spectrum of $A$ outside the closed unit disk:
  \begin{equation}
  \label{projector}
    P=I_N-\frac{1}{2\pi i}\int_{|z|=1}(zI_N-A)^{-1}dz.
  \end{equation}
  Then $f(z)=\sum_{n\in\mathbb Z}f_nz^n$ with
  \begin{equation}
    f_n=\begin{cases}\, CA^{n-1}(I_N-P)B,\quad\, n>0,\\
      \,D\delta_{n0}-CA^{n-1}PB,\quad n\le 0,
      \end{cases}
  \end{equation}
  and $\sum_{\mathbb Z} \|f_n\|<\infty$.
  \label{Arche-de-la-defense}
\end{theorem}

We now turn to the $\sharp$ symmetric case. Let us consider $J\in\mathbb{C}^{2\times 2}$ as above, with $J\, J^*=I_2$. 
\begin{theorem}
\label{sharp_a}
  Let $f(z)$ be a $\mathbb C^{2p\times 2p}$-valued rational function, analytic at infinity, $\sharp$-symmetric with respect to $J$, and regular on the unit circle.
  Then $f$ admits a $\sharp$-symmetric realization, analytic at infinity, and regular on the unit circle.
\end{theorem}

\begin{proof}
Going back to the construction in the proof of Theorem \ref{real} we see that the matrix coefficients $H_{k,m}$ are $\sharp$-symmetric, and so is the corresponding realization as written in Theorem~\ref{real}. 
We note that the Riesz projection $P$ obtained in Theorem~\ref{Riesz_proj} is also $\sharp$-symmetric.
        \end{proof}

        \begin{theorem}
        \label{sharp_b}
          Let $f\in\mathcal W^{2p\times 2p}$ rational and taking strictly positive values and $\sharp$ symmetric. Then it admits rational
          $\sharp$-symmetric left and right factorizations.
        \end{theorem}

        \begin{proof}
This follows from the fact that the a $\sharp$-symmetric rational function admits a sharp-symmetric realization. 
          \end{proof}

Assuming $f(e^{it})>0$, one can write $f(z)=w(z)(w(1/\overline{z}))^*$ where the spectral factor $w(z)$ is rational.
For the following,  see also \cite[Lemma 1.2, p. 145]{MR1294714} where the realization is centered at the origin. The computations appear also in
\cite{ACDSV}.
\begin{proposition}
  \label{bruxelles}
Let 
\begin{equation*}
d+c(z I_m-a)^{-1}b
\end{equation*}
be  a realization of  the  left spectral factor of $f(z)$ with both $\sigma(a)$ and $\sigma(a-bd^{-1}c)$ inside the open unit disk.
Then,
\begin{equation*}
  f(z)=D+C(zI_{2m}-A)^{-1}B,
\end{equation*}
where
\begin{align*}
A&=\begin{pmatrix}a&-bb^*a^{-*}\\0&a^{-*}\end{pmatrix}\\
B&=\begin{pmatrix}b(d^*-b^*a^{-*}c^*)\\ a^{-*}c^*\end{pmatrix}\\
   C&=\begin{pmatrix}c&-db^*a^{-*}\end{pmatrix}\\
  D&=d(d^*-b^*a^{-*}c^*).
\end{align*}
\end{proposition}

The Riesz projection $P$ from~\eqref {projector} becomes:
\begin{equation}
\begin{pmatrix}I_m&X\\0&0\end{pmatrix}
\end{equation}
where $X$ is the unique solution to the Stein equation
\begin{equation}
  X-aXa^*=bb^*
\end{equation}
We now give formulas for the Fourier coefficients in terms of $a,b,c$ and $d$. 

\begin{proposition}
    \label{r0}
  The Fourier coefficients   of the spectral function  are given by
  \begin{align*}
  f_0&=dd^*+cXc^*\\
  f_k&=(db^*+cXa^*)a^{*(k-1)}c^*,\quad k=1,2,\ldots\\
  f_{-k}&=f_k^*,\quad k=1,2,\ldots
  \end{align*}
  in terms of a realization \eqref{bruxelles} of the spectral factor.
\end{proposition}

Using Definition~\ref{bc_M_sharp}, one can re-write Proposition~\ref{bruxelles} and Proposition~\ref{r0} in terms of the $\sharp$ conjugate as well.
\
\begin{remark}
Let us assume that $f(e^{it})>0$ is also $\sharp$-symmetric and the dimensions are even. It follows that the spectral factor $w(z)$, which yields the decomposition of $f$ i.e. $f(z)=w(z)(w(1/\overline{z}))^*$, is rational and $\sharp$ symmetric.
In this context, in Propositions~\ref{bruxelles} and~\ref{r0} all factors $a, b, c, d$ as well as $A, B, C, D$ will be $\sharp$ symmetric.
The Riesz Projector $P$ will be $\sharp$-symmetric as well.
  \end{remark}      
        
The results in this section will be interpreted in Section 4 in the setting of the $\mathbb{BC}$ Wiener algebra.

\subsection{Wiener algebras and classical Superoscillations}
\label{classical_super}
The notion of superoscillations originate with the works of Aharonov and Berry and the notion of weak measurements. We refer to the papers
\cite{aav,av,abook,b1,b4,MR4:218g}. We refer to \cite{AOKI,QS2,distributions} for
recent developments and to the paper \cite{acs_jfaa} for recent connections with Schur analysis.



\begin{definition}
A superoscillatory sequence is a sequence of complex-valued functions $F_m(t,a)$ defined on $\mathbb R$ as follows:
\begin{equation}
\label{sup_F_n}
F_m(t,a)=\left(\cos(\frac{t}{m})+i\, a \cos(\frac{t}{m})\right)^m=\sum_{k=0}^m c_k(m,a) e^{it(1-2k/m)},
\end{equation}
where:
\begin{equation}
\label{sup_F_coefs}
c_m(n,a)=\binom{m}{k}\left(\frac{1+a}{2}\right)^{m-k}\left(\frac{1-a}{2}\right)^{k},
\end{equation}
with $a>1$.
\end{definition}

One can then see that for fixed $t\in\mathbb R$ we have:
\begin{equation}
\label{sup_seq}
\lim_{m\to\infty} F_m(t,a)=e^{iat}
\end{equation}
and convergence is uniform on compact subsets of the real line. This sequence plays an important role in the respective Wiener algebra as follows.

\

Let us now turn to the matrix-valued Wiener algebra $\mathcal W^{p\times p}$ case. We will show that the superoscillatory sequence $F_m(t,a)$ provides a good approximation in this case as well.
We have the following result which will expand to the $\mathbb{BC}$ as well.

\begin{theorem}
\label{Suposc_approx}
  Let $f\in\mathcal W^{p\times p}$, then $f$ can be uniformly approximated on compact subsets of $\mathbb{R}$ by
  $$
  f_{-1}e^{-it}+f_0+ f_1e^{it}+\sum_{\substack{n\in\mathbb Z\\ n\not\in\left\{-1,0,1\right\}}} f_nF_m(t,n).
  $$ 
  
\end{theorem}

\begin{proof}
Let $\epsilon>0$ and $f\in\mathcal W^{p\times p}$. If we write
\[
f(e^{it})=f_0+f_{-1}e^{-it}+f_1e^{it}+g(e^{it}),
\]
there exists $N\in\mathbb N$ such that
\[
\sum_{|n|>N}\|f_n\|<\frac{\epsilon}{2},
\]
and by~\eqref{sup_seq} there exists $M$ such that
\[
  \forall n\in\left\{-N,\ldots,N\right\}\setminus\left\{-1,0,1\right\},\quad m\ge M\,\,\Longrightarrow\,\, |F_m(t,n)-e^{int}|
  \le \frac{\epsilon}{2N-2}.
\]
Thus for $m\ge M$
\[
\|g(e^{it})-\sum_{\substack{n\in\mathbb Z\\ n\not\in\left\{-1,0,1\right\}}} f_nF_m(t,n)\|\le \epsilon,
\]
and so
\[
\left\|f(e^{it})-\left(f_{-1}e^{-it}+f_0+f_1e^{it}+\sum_{\substack{n\in\mathbb Z\\ n\not\in\left\{-1,0,1\right\}}} f_nF_m(t,n)\right)\right\|\le \epsilon
\]  
\end{proof}

We will extend these results to the bicomplex case using the intrinsic structure of its algebra and analysis.


\section{Bicomplex Analysis and Symmetry Domains}
\label{bcplx}

In this section we revisit the analysis on the bicomplex algebra and write the appropriate symmetry domains that will allow us to introduce a
bicomplex Wiener algebra in this case.

\subsection{Bicomplex algebra as $2\times 2$ complex matrices}
\label{bc-const}
As described in~\cite{ALV}, in a way similar to the space of complex numbers and the Pauli model, one can re-write the space of bicomplex numbers as a subspace of complex matrices $\mathbb C^{2\times 2}$ as:
\begin{equation}
  \label{bcform}
  \mathbb{BC}=\left\{  \begin{pmatrix}z_1&-z_2\\ z_2&z_1\end{pmatrix}\,;\, z_1,z_2\in\mathbb C \right\},
  \end{equation}
  where the complex unit in $\mathbb C$  is denoted by $i$.
We write
  \[
    \begin{pmatrix}z_1&-z_2\\ z_2&z_1\end{pmatrix}=z_1\begin{pmatrix}1&0\\ 0&1\end{pmatrix}+z_2\begin{pmatrix}0&-1\\ 1&0\end{pmatrix}=z_1\begin{pmatrix}1&0\\ 0&1\end{pmatrix}+J \,z_2,
  \]
  with $J$ as in~\eqref{bc-J}.
  We will use the shorter notation $Z=z_1+jz_2$ to denote a bicomplex number. For example, in this writing, we have that the original
  complex unit $i$ becomes:
   \[
    i\,I=\begin{pmatrix} i&0\\ 0&i\end{pmatrix}.
  \]

\begin{remark}
 In the same notation as in Subsection~\ref{bcplx_alg_intro} one obtains the corresponding matrix to the {\em hyperbolic unit} $k=i\,j,$ (i.e. $k^2=1$) as $K=i J$.
 \end{remark}
 
 In the matrix notation $k$ is represented by:
 
  \[
  i\,J=\begin{pmatrix} 0&-i\\ i&0\end{pmatrix},
  \]
and one can easily check that the square of this matrix is the identity matrix.\smallskip


In the matricial writing of a bicomplex number $Z\in\BC$, we make use of the following primary decomposition and the fact that a normal matrix is unitary diagonalizable to obtain:
\[
  \begin{pmatrix}z_1&-z_2\\ z_2&z_1\end{pmatrix}=\frac{1}{\sqrt{2}} \begin{pmatrix}1&1\\ i&-i\end{pmatrix} \begin{pmatrix}z_1-iz_2&0\\ 0&z_1+iz_2\end{pmatrix}\frac{1}{\sqrt{2}}
      \begin{pmatrix}1&-i\\1&i\end{pmatrix}.
    \]
    The unitary diagonalization matrix is:
    \begin{equation}
      \label{U-diag}
U=\frac{1}{\sqrt{2}} \begin{pmatrix}1&1\\ i&-i\end{pmatrix},
\end{equation}
and the eigenvalues are:
\begin{equation}
\label{lambdas}
\begin{split}
  \lambda_1&=z_1-iz_2\\
  \lambda_2&=z_1+iz_2.
\end{split}
\end{equation}

\begin{remark}
In this commutative setting we have the same zero divisors $\e_1,\,\e_2$ as in Subsection~\ref{bcplx_alg_intro}, given by the diagonalization.
As expected, these split the bicomplex space in $\BC=\mathbb C \mathbf{e}_1\bigoplus \mathbb C \mathbf{e}_2$, as:
\begin{equation}
  Z=z_1+jz_2=(z_1-iz_2)\mathbf{e}_1+(z_1+iz_2)\mathbf{e}_2.
\end{equation}
\end{remark}

Using normality again, the corresponding matrix form of a bicomplex number can be written as the associated sum of two
weighted orthogonal projections
\[
\begin{pmatrix}
z_1&-z_2\\z_2&z_1
\end{pmatrix}
=
\lambda_1\cdot E_1+\lambda_2\cdot E_2
\]
where the eigenvalues $\lambda_{1,2}$ are given in~\eqref{lambdas}  and the orthogonal projections are
\[
E_1=\frac{1}{2}\begin{pmatrix}1&-i\\i&1\end{pmatrix}\quad
E_2=\frac{1}{2}\begin{pmatrix}1&i\\-i&1\end{pmatrix}.
\] 

It is easily seen that $E_1$ and $E_2$ are the corresponding matrix forms of $ \mathbf{e}_1$ and $ \mathbf{e}_2$.

\

We will now define the conjugates in the bicomplex setting, as in~\cite{CSVV,bcbook}.

\begin{definition} For any $Z\in \BC$ we have the following three conjugates:
  \begin{eqnarray}
  \label{conj}
    \overline{Z}=\overline{z_1}+j\overline{z_2}\\
     Z^{\dagger}=z_1-jz_2\\
      Z^*=\overline{Z^{\dagger}}=\overline{z_1}-j\overline{z_2}.
  \end{eqnarray}
\end{definition}

\begin{remark}
In the matrix form the space of {\em hyperbolic numbers $\mathbb D$} is realized as the subset of matrices \eqref{bcform} with:
  \[
\begin{pmatrix}a&-ib\\ ib&a\end{pmatrix}.
  \]  
   \end{remark}


As in~\cite{ALSS}, one can define the corresponding notions of positivity:

\begin{definition}
\label{BC_pos_idem}
  A matrix $M\in\mathbb{BC}$ is called {\em Hermitian} if $M=M^*$ and {\em positive} if $c^*Mc\ge 0$ for all $c\in\mathbb{BC}^n$.
  Here $M^*$
\end{definition}

In~\cite{ALV} we proved  the following equivalent definition, in terms of our matricial form:
\begin{proposition}
\label{ALV_hermitian_positive}
  Let $M=M_1+jM_2\in\mathbb{BC}^{n\times n}$ with $M_1,M_2\in\mathbb C^{n\times n}$. Then $M$ is Hermitian, (resp. positive), in the sense of Definition \ref{BC_pos_idem} if and only if
\[
  \begin{pmatrix}M_1&-M_2\\ M_2&M_1\end{pmatrix}
\]
is Hermitian (resp. a positive matrix, in the sense of positive matrices with complex entries).
If we write
  \[
M_1+jM_2=P_1\mathbf{e}_1+P_2\mathbf{e}_2,
\]
with $P_1=M_1-iM_2$ and $P_2=M_1+iM_2$.
Then, $M$ is positive in the sense of Definition \ref{BC_pos_idem} if and only if both $P_1$ and
$P_2$ are positive  elements of $\mathbb C^{n\times n}$. 
\end{proposition}

\subsection{$\mathbb{BC}$ analyticity:}

As seen in~\cite{CSVV,bcbook} and~\cite{mltcplx}, the function $F=F_1+jF_2$ is $\mathbb{BC}-$analytic iff and only if the functions $F_1$ and $F_2$ are complex holomorphic in $z_1$ and $z_2$ i.e.:
\begin{equation}
\frac{\partial F_1}{\partial \overline{z_1}}=\frac{\partial F_1}{\partial \overline{z_2}}=\frac{\partial F_2}{\partial \overline{z_1}}=\frac{  \partial F_2}{\partial \overline{z_2}}=0,
\end{equation}
and if the following Cauchy-Riemann like equations hold:
\begin{eqnarray}
  \label{cr1}
\frac{  \partial F_1}{\partial z_1}&=&\frac{  \partial F_2}{\partial z_2}\\
  \frac{\partial F_1}{\partial z_2}&=&-\frac{  \partial F_2}{\partial z_1}.
                                         \label{cr2}
\end{eqnarray}
\smallskip
\begin{remark}
The identity function $Z=z_1+jz_2$ is $\mathbb{BC}-$analytic, but the $\mathbb C$-valued functions $Z\mapsto z_1$ and $Z\mapsto z_2$ are not.
\end{remark}
\smallskip


\subsection{Bicomplex Matrices as Complex Matrices with symmetry}
\label{bcplx_matrix_sharp}

 \  
We will now view bicomplex matrices in $\BC^{p\times p}$ as a matrices in $\C^{2p\times 2p}$ with the bicomplex $\sharp-$symmetry defined below.
We then can reduce the various results to ones in a complex Wiener algebra with symmetry.
Note that we will also use the second approach, splitting the problem into two classical complex ones via the idempotent representation. 

\begin{remark}{\rm     
This type of symmetry can be viewed as a complexification of the geometry of complex numbers as well. While in the complex framework a rotation matrix takes the Pauli form, in the bicomplex case it becomes a "double rotation".
}
\end{remark}

We recall that in the case the bicomplex algebra, as seen is~\cite{ALV} and in~\eqref{bc-J}, we have an important structural matrix :
$$
J=\begin{pmatrix} 0 &-1\\ 1&0\end{pmatrix}\, ,
$$
which has the obvious property that $J\, J^*=I_2$.
In this context, from Definition~\ref{bcplx-sharp-2by2}, we have:
\begin{definition}
The matrix $M\in\mathbb C^{2\times 2}$ is called bicomplex $\sharp$-symmetric if it satisfies $M=M^\sharp$, with respect to this $J$.
\end{definition}

This matrix gives us the structure of bicomplex numbers as $M$ is $J$-symmetric if and only if it is of the form (see~\ref{petunia}).
\begin{equation*}
M=\begin{pmatrix}z_1&-z_2\\ z_2&z_1\end{pmatrix}\,;\, z_1,z_2\in\mathbb C.
\end{equation*}

For the general matrix case, defining $$J_a=J \otimes I_a\,=\begin{pmatrix} 0 &-1\\ 1&0\end{pmatrix}\, \otimes I_a,$$
and $J_a^*=J^* \otimes I_a$, where $\otimes$ denotes the  Kronecker (tensor) product of matrices, we have:

\begin{definition}
\label{bcplx-sharp-2p_by_2p}
The block matrix $M\in\mathbb C^{2a\times 2b}$ is called bicomplex $\sharp$-symmetric if and only if it satisfies  $M=M^\sharp$, where
\begin{equation*}
M^\sharp=J_aMJ_b^*.
\end{equation*}
\end{definition}

We will make use of all of these notions in the following Section~\ref{Wiener_BC_matrix}, specifically in Subsection~\ref{bcplx_wiener_sharp}.
\
The following lemma is an immediate result:
\begin{lemma}
\label{pos_two_structures}
Let $M\in \BC^{p\times p}$, with $M=(Z_{st})_{1\leq s,t\leq p}$ and $Z_{st}=z^1_{st}+j z^2_{st}\in\BC$. 
We can decompose $M$ as a matrix in $\C^{2p\times 2p}$ first as in Proposition~\ref{ALV_hermitian_positive}
and second as a matrix where each $Z_{st}=z^1_{st}+j z^2_{st}$ is written as the $2\times 2$ block 
\[
  \begin{pmatrix}z^1_{st}&-z^2_{st}\\ z^2_{st} &z^1_{st}\end{pmatrix}.
\]
Then the second matrix is bicomplex $\sharp$ symmetric and the definition of positivity coincide.
\end{lemma}

\
We now turn our attention to building the various bicomplex Wiener algebras.

\section{The Scalar and Matrix-Valued Bicomplex Wiener Algebra}
\label{Wiener_BC_matrix}
\setcounter{equation}{0}

We extend the Wiener algebra setting to the set of bicomplex numbers in two ways, as follows.
The bicomplex Wiener algebra can be defined in a direct way using the idempotent split, or as a symmetry based Wiener algebra in the complex domain, using a block matrix that fixes the bicomplex structure.

\subsection{Scalar Bicomplex Wiener Algebras}
\label{scalar_W_BC}
Here are the theorems referenced in  the introduction that show the structure of the scalar bicomplex valued case:
\begin{theorem}
Assume that $f(Z)$ invertible on $\partial\mathbb K$. Then $f$ is invertible in the bicomplex Wiener algebra $\mathfrak W$.
\end{theorem}

\begin{theorem}
Assume that $f(Z)>0$ invertible on $\partial\mathbb K$. Then $f$ admits left and right spectral factorizations in $\mathfrak W_+$ and $\mathfrak W_-$. 
  \end{theorem}
The proofs of these theorems follow directly from Theorem~\ref{p-bcplx-wiener}  and Theorem~\ref{wiener_m_sym_complex}, in the case when $p=1$.

\subsection{Wiener bicomplex algebra using bicomplex analyticity}
\setcounter{equation}{0}
\label{sec-4-2}

\
In this section we describe the general matrix-valued bicomplex Wiener algebra, induced by the analytic and algebraic structures of the set of bicomplex numbers, via the idempotent split.

We first recall the definition of matrix-valued bicomplex Wiener algebras, see Definition~\ref{bc_wiener-matrix}.
 \begin{definition}
The bicomplex Wiener algebra $\mathfrak W^{p\times p}$ consists of the matrix-valued functions of the form 
\begin{equation}
  f(Z)=\sum_{n\in\mathbb Z} f_nZ^n
\end{equation}
with $Z\in \partial\mathbb K$, $f_n\in\BC^{p\times p}$ and
$\sum_{n\in\mathbb Z}||f_n||<\infty$.

Here $||f_n||$ denotes the following norm of the bicomplex matrix $f_n=f_{n,1}\e_1+f_{n,2}\e_2$:
$||f_n||=\|f_{n,1}\|+\|f_{n,2}\|$, where $\|f_{n,1}\|$ and $\|f_{n,2}\|$ are the operator norms as complex matrices.

This Wiener algebra is endowed with pointwise multiplication and norm
\begin{equation}
  \|f\|=\sum_{n\in\mathbb Z}\|f_n\|.
\end{equation}
\end{definition}

The elements of $\mathfrak W^{p\times p}_+$  are  of the form 
\begin{equation}
  f(Z)=\sum_{n\ge 0} f_nZ^n,
\end{equation}
while elements of $\mathfrak W^{p\times p}_-$  are of the form
 \begin{equation}
  f(Z)=\sum_{n\ge 0} f_{-n}Z^{-n}.
\end{equation}

We have the following theorems that describe the structure of the Wiener algebras in this setting.

\begin{theorem}
\label{p-bcplx-wiener}
Let $f\in \mathfrak W^{p\times p}$ be such that $f(Z)$ is invertible for $Z\in\partial \mathbb K$. Then $f$ is invertible in
  $\mathfrak W^{p\times p}$.
\end{theorem}

\begin{proof}
We write $f(Z)=\sum_{n\in\mathbb Z} f_nZ^n$ with $f_n\in\mathbb{BC}^{p\times p}$, and  $f_n=f_{1n}e_1+f_{2n}e_2$ where
$f_{1n}$ and $f_{2n}$ belong to $\mathbb C^{p\times p}$.
With $Z\in\partial\mathbb K$ given by \eqref{z-b-a-1} and corresponding $\lambda_1$ and $\lambda_2$ given by \eqref{dupont-circle-1}
  \[
\sum_{n\in\mathbb Z}f_{1 n}e^{in(b-a)}>0\quad {\rm and}\quad \sum_{n\in\mathbb Z}f_{2 n}e^{in(b+a)}>0,\quad a,b\in[0,2\pi].
\]
Setting $a=0$, the classical inversion result leads to
\[
\left(\sum_{n\in\mathbb Z}f_{\ell n}e^{inb}\right)^{-1}=g_\ell(e^{ib})\in\mathcal W^{p\times p},\quad \ell=1,2.
\]
We set
\[
g(Z)=\left(\sum_{n\in\mathbb Z} g_{1 n}\lambda_1^n\right) {\mathbf e_1}+\left(\sum_{n\in\mathbb Z}g_{2 n}\lambda_2^n \right){\mathbf e_2}.
  \]
  Then $g\in(\mathcal W(\mathbb{BC}))^{p\times p}$ and $f(Z)g(Z)=I_p$ for $Z\in\partial \mathbb K$.
\end{proof}

\begin{corollary}
In the previous theorem it is enough to require invertibility on $e^{it}$, $t\in[0,2\pi]$ to insure invertibility in $\mathfrak W^{p\times p}$.
\end{corollary}

\begin{theorem}
\label{bc_m_sym_complex}
  Let $f\in \mathfrak{W}^{p\times p}$ be such that $f(Z)>0$ for $Z\in\partial \mathbb K$. Then
  $f$ admits a spectral factorization $f(Z)=f_+(Z)f_+(Z)^*$ where $f_+^{\pm1}\in \mathfrak W^{p\times p}_+$.
\end{theorem}

\begin{proof}
  In the same notation as in the proof of the previous theorem, and with $Z\in\partial\mathbb K$ given by \eqref{z-b-a} and corresponding
$\lambda_1$ and $\lambda_2$ given by \eqref{dupont-circle}
  \[
\sum_{n\in\mathbb Z}f_{1 n}e^{in(b-a)}>0\quad {\rm and}\quad \sum_{n\in\mathbb Z}f_{2 n}e^{in(b+a)}>0,\quad a,b\in[0,2\pi].
\]
Setting $a=0$, the classical factorization result leads to the spectral factorization
\[
\sum_{n\in\mathbb Z}f_{\ell n}e^{inb}=f_{\ell +}(e^{ib})(f_{\ell +}(e^{ib}))^*,\quad \ell=1,2.
\]
with $f_{\ell +}\in\mathcal W_+^{p\times p}$ and invertible in $\mathcal W_+^{p\times p}$. Let
\[
  f_{\ell +}(e^{ib})=\sum_{n=0}^\infty g_{\ell, n}e^{inb},\quad \ell=1,2,
\]
with $\sum_{n=0}^\infty\|g_{\ell,n}\|<\infty$ for $\ell=1,2$, and set
\[
f_+(Z)=\left(\sum_{n=0}^\infty g_{1n}e^{in(b-a)}\right)\mathbf e_1+\left(\sum_{n=0}^\infty g_{2n}e^{in(b+a)}\right)\mathbf e_2, \quad Z=e^{ib+ja}\in\partial\mathbb K.
\]
Then, $f(Z)=f_+(Z)f_+(Z)^*$. To conclude we need to verify that $f_+$is invertible in $\mathfrak W^{p\times p}_+$. Let
\[
(f_{\ell +}(e^{ib}))^{-1}=\sum_{n=0}^\infty h_{\ell, n}e^{inb},\quad \ell=1,2,
\]
with $\sum_{n=0}^\infty|h_{\ell,n}|<\infty$ for $\ell=1,2$. We have
\[
(f_+(Z))^{-1}=\left(\sum_{n=0}^\infty h_{1n}e^{in(b-a)}\right)\mathbf e_1+\left(\sum_{n=0}^\infty h_{2n}e^{in(b+a)}\right)\mathbf e_2, \quad Z=e^{ib+ja}\in\partial\mathbb K,
\]
so that $f_+(Z)$ is invertible in $\mathfrak W^{p\times p}_+$.
\end{proof}

\begin{corollary}
In the previous theorem it is enough to assume that $f(e^{ib})>0$ for $b\in[0,2\pi]$.
\end{corollary}


\subsection{Bicomplex Wiener Algebras using $\sharp$ symmetries}
\label{bcplx_wiener_sharp}
Let us return to bicomplex matrices with $\sharp$ symmetry as in Definition~\ref{bcplx-sharp-2p_by_2p}.

  \begin{proposition}
  \label{bcplx_sharp_prop}
Using our specific bicomplex $\sharp$ definition (with $J$ defined in~\eqref{bc-J}) and setting $p,q,s\in\mathbb N$ and $f\in\mathcal W^{2p\times 2q}$ and $g\in\mathcal W^{2q\times 2s}$, we have:
    \begin{eqnarray}
      (fg)^\sharp(e^{it})&=&f^\sharp(e^{it})g^\sharp(e^{it})\\
      f^\sharp(e^{it})&=&(f(e^{it}))^\sharp\\
      ((f(e^{it}))^\sharp)^*&=&((f(e^{it}))^*)^\sharp.
    \end{eqnarray}
    Furthermore, 
$f\in\mathcal W_+^{2p\times 2q}$ if and only if $f^\sharp\in\mathcal W_+^{2p\times 2q}$.
    \end{proposition}

\begin{proposition}
  The function $f\in\mathcal  W^{2p\times 2p}$ is bicomplex $\sharp$-symmetric if and only if its Fourier coefficients are bicomplex $\sharp$-symmetric:
  \[
f_n=(f_n)^\sharp,\quad n\in\mathbb Z.
    \]
  \end{proposition}

\begin{theorem}
\label{wiener_m_sym_bicomplex-2}
  Let $f\in\mathcal W^{2p\times 2p}$ be bicomplex $\sharp$-symmetric and such that $f(e^{it})>0$ for $t\in[0,2\pi]$. Then, there is a uniquely determined spectral
  factorization $f=f_+f_+^*$ where $f_+\in\mathcal W_+^{2p\times 2p}$ is bicomplex $\sharp$-symmetric and satisfies $f_+(1)=I_{2p}$ and invertible and any  other spectral factorization differs by a bicomplex $\sharp$-symmetric
  unitary multiplicative constant on the right from $f_+$.
\end{theorem}

Using Lemma~\ref{pos_two_structures}, we can now establish the equivalence between the two types of Wiener algebras in the bicomplex case.
\begin{theorem}
\label{W_BC_equiv}
$\mathfrak W^{p \times p}$ is isomorphic to the space of $\mathcal W^{2p \times 2p}$ with bicomplex $\sharp$ symmetry.
\end{theorem}

\subsection{Rational Bicomplex functions and the Bicomplex Wiener Algebra}
\label{rational_bcplx_sharp}
\setcounter{equation}{0}


We now re-write Theorem~\ref{sharp_a} and Theorem~\ref{sharp_b} in the bicomplex case. The proofs are left to the reader.

\begin{theorem}
  Let $f(z)$ be a $\mathbb C^{2p\times 2p}$-valued rational function, analytic at infinity, bicomplex $\sharp$-symmetric, and regular on the unit circle.
  Then $f$ admits a bicomplex $\sharp$-symmetric realization, analytic at infinity, and regular on the unit circle.
\end{theorem}

\begin{proof}
Follows from proof of Theorem~\ref{sharp_a}.

      \end{proof}

        \begin{theorem}
          Let $f\in\mathcal W^{2p\times 2p}$ rational and taking strictly positive values and bicomplex $\sharp$ symmetric. Then it admits rational
         bicomplex  $\sharp$-symmetric left and right factorizations.
        \end{theorem}

        \begin{proof}
This follows from the fact that the a bicpomplex $\sharp$-symmetric rational function admits a bicomplex $\sharp$-symmetric realization. 
          \end{proof}

\section{Bicomplex Superoscillations}\
\label{Bcplx_sup}

As in Lemma~\ref{bcpx_trig}, and formulas~\ref{ab-ab} and~\ref{z-b-a}, the distinguished boundary $|\lambda_1|=|\lambda_2|=1$ of the bidisk in the decomposition along the idempotents corresponds to
  \begin{equation*}
    \left\{ e^{it}(\cos s, \sin s)\, ;t,s\in[0,2\pi]\right\}, 
  \end{equation*}
  i.e.
  \begin{equation}
    Z=e^{it}e^{js}
    \end{equation}
and we have for $z_1=uz_2$ when $z_2\not=0$:
\[
\partial \mathbb K=  \left\{(e^{i\theta},0), \theta\in[0,2\pi]\right\}\bigcup\left\{\bigcup_{u\in\mathbb R}\left\{\left(\frac{ue^{i\theta}}{\sqrt{u^2+1}},\frac{e^{i\theta}}{\sqrt{u^2+1}}\right),\,
    \theta\in[0,2\pi]\right\}\right\},
\]
which is rewritten as \eqref{ab-ab}. To conclude, note that
\[
z_1+jz_2=e^{it}(\cos s +j\sin s)=e^{it}e^{js},
\]  

and

  {\rm
    \begin{equation}
      \label{dupont-circle}
    \lambda_1=e^{i(t-s)}\quad{\rm and}\quad     \lambda_2=e^{i(t+s)}
  \end{equation}
  and
  \[
    z_1+jz_2=e^{i(t-s)}e_1+e^{i(t+s)}e_2
  \]
  }

In the bicomplex case, for given $x,y,a,b\in\R$ and $Z=\lambda_1\e_1+\lambda_2\e_2$, where $$\lambda_1=\cos(\frac{x}{m})+i \, a \cos(\frac{x}{m})$$ and $$\lambda_2=\cos(\frac{y}{m})+i \, b \cos(\frac{y}{m})$$ and we can write the following bicomplex superoscillatory sequence:

\begin{eqnarray}
\label{bc_sup_F_m}
F_m(x,y,a,b)=&(\lambda_1 \e_1+\lambda_2 \e_2)^m=\lambda_1^m\e_1+\lambda_2^m\e_2\\
=&\left(\cos(\frac{x}{m})+i \, a \cos(\frac{x}{m})\right)^m\e_1+\left(\cos(\frac{y}{m})+i \, b \cos(\frac{y}{m})\right)^m\e_2.
\end{eqnarray}

Using the respective Lemmas in the complex case we have the following structure of the bicomplex superpscillations:

\begin{lemma}
We have that:
\begin{equation*}
F_m(x,y,a,b)=\sum_{k=0}^n c_k(m,a) e^{ix(1-2k/m)}\e_1+c_k(m,b) e^{iy(1-2k/m)}\e_2,
\end{equation*}
where 
\begin{equation}
\label{bcsup_coeffs}
c_k(m,a)=\binom{m}{k}\left(\frac{1+a}{2}\right)^{m-k}\left(\frac{1-a}{2}\right)^{k},\quad
c_k(m,b)=\binom{m}{k}\left(\frac{1+b}{2}\right)^{m-k}\left(\frac{1-b}{2}\right)^{k}.
\end{equation}
\end{lemma}

We also have the same type of limiting behaviour:
\begin{lemma}
We have that:
$$
\lim_{m\to\infty} F_m(x,y,a,b)=e^{iax}\e_1+e^{iby}\e_2,
$$
and $F_m(x,y,a,b)=F_m(x,a)\e_1+F_m(y,b)\e_2$.
\end{lemma}

\begin{theorem}
  Let $f\in\mathfrak W^{p\times p}$, then $f(Z)=\sum_{n\in\mathbb Z} f_nZ^n,$ where $Z=e^{i t}e^{j s}$ can be uniformly approximated on compact subsets of $\mathbb{R}$ by ....
  
\end{theorem}

\begin{proof}
We have that 
$$f(Z)=\sum_{n\in\mathbb Z} f_n(e^{i t}e^{j s})^n=\sum_{n\in\mathbb Z} f_n(e^{i (t-s)}\e_1+e^{i (t+s)}\e_2)^n$$
Therefore

$$
f(Z)=\sum_{n\in\mathbb Z} (A_n\e_1+B_n\e_2)(e^{i (t-s)}\e_1+e^{i (t+s)}\e_2)^n
$$

$$
f(Z)=\sum_{n\in\mathbb Z} A_n e^{i n(t-s)}\e_1+B_ne^{i n(t+s)}\e_2=: f^1(e^{i(t-s)})\e_1+f^2(e^{i(t+s)})\e_1
$$

Using Theorem~\ref{Suposc_approx} we have that $f_1(e^{i(t-s)})$ is approximated by: 

$$\mathcal E_1= f^1_{-1}e^{-i(t-s)}+f^1_0+ f^1_1e^{i(t-s)}+\sum_{\substack{n\in\mathbb Z\\ n\not\in\left\{-1,0,1\right\}}} f^1_nF_m(t-s,n).$$

Similarly, $f_2(e^{i(t+s)})$ is approximated by:

$$\mathcal E_2= f^2_{-1}e^{-i(t+s)}+f^2_0+ f^2_1e^{i(t+s)}+\sum_{\substack{n\in\mathbb Z\\ n\not\in\left\{-1,0,1\right\}}} f^2_nF_m(t+s,n).$$

Using the Lie norm induced on the space of bicomplex matrices we can write an approximation for $f(Z)$ to be $\mathcal E_1\e_1+ \mathcal E_2\e_2$:

\begin{eqnarray*}
\displaystyle
f(Z)= &(f^1_{-1}e^{-i(t-s)}+f^1_0+ f^1_1e^{i(t-s)})\e_1+(f^2_{-1}e^{-i(t+s)}+f^2_0+ f^2_1e^{i(t+s)})\e_2 +\\
 	&+\sum_{\substack{n\in\mathbb Z\\ n\not\in\left\{-1,0,1\right\}}} f^1_nF_m(t-s,n) \e_1+ \sum_{\substack{n\in\mathbb Z\\ n\not\in\left\{-1,0,1\right\}}} f^2_nF_m(t+s,n) \e_2.
\end{eqnarray*}
which yields:
\begin{equation*}
f(Z)= f_{-1}e^{-it}e^{-js}+f_0+ f_1e^{it}e^{js}+ \sum_{\substack{n\in\mathbb Z\\ n\not\in\left\{-1,0,1\right\}}} f_n(F_m(t-s,n)\e_1+F_m(t+s,n )\e_2),
\end{equation*} 
namely:
\begin{equation*}
f(Z)= f_{-1}e^{-it}e^{-js}+f_0+ f_1e^{it}e^{js}+ \sum_{\substack{n\in\mathbb Z\\ n\not\in\left\{-1,0,1\right\}}} f_n\,F_m(t-s,t+s,n,n).
\end{equation*} 

This completes the proof.
\end{proof}

\
This application concludes the paper and we expect to be able to apply this framework to other settings as well. Wiener Algebras have many applications in the theory of Linear Systems in the classical case, so we expect results in this direction as well. 
 
 We are also investigating whether the setting of bicomplex Wiener algebras is useful in other applications of Quantum Mechanics.


\end{document}